\documentclass[12pt,eqno,]{article}
\usepackage{amsmath,amsthm,amssymb}

\makeatletter

\newcommand{\Rmnum}[1]{\expandafter\@slowromancap\romannumeral #1@}
\makeatother

\newcommand{\arcsinh}{\rm {arcsinh}}
\title{\normalsize\bf  SHARP INEQUALITIES FOR THE NEUMAN-S\'{A}NDOR MEAN IN TERMS OF ARITHMETIC AND CONTRA-HARMONIC MEANS}
\author{\small  MIAO-KUN WANG, YU-MING CHU AND BAO-YU LIU}
\date{}

\begin{document}
\maketitle

\renewcommand{\thefootnote}{\fnsymbol{footnote}}

\footnotetext{\hspace*{-5mm}
\begin{tabular}{@{}r@{}p{14.0cm}@{}}
&\qquad Mathematics Subject Classification (2010): 26E20.\\
&\qquad Keywords and phrases: Neuman-S\'{a}nder mean, arithmetic mean, contra-harmonic mean.\\
&\qquad This research was supported by the Natural Science
Foundation of China under Grants 11071069 and 11171307, and Innovation
Team Foundation of the Department of Education of Zhejiang Province
under Grant T200924.
\end{tabular}}

\vspace{-1cm}

\def\abstractname{}
\begin{abstract}
\noindent {\emph{Abstract}.} In this paper, we find the greatest
values $\alpha$ and $\lambda$, and the least values $\beta$ and
$\mu$ such that the double inequalities
$$C^{\alpha}(a,b)A^{1-\alpha}(a,b)<M(a,b)<C^{\beta}(a,b)A^{1-\beta}(a,b)$$
and
\begin{align*}
&[C(a,b)/6+5 A(a,b)/6]^{\lambda
}\left[C^{1/6}(a,b)A^{5/6}(a,b)\right]^{1-\lambda}<M(a,b)\\
&\qquad<[C(a,b)/6+5
A(a,b)/6]^{\mu}\left[C^{1/6}(a,b)A^{5/6}(a,b)\right]^{1-\mu}
\end{align*}
hold for all $a,b>0$ with $a\neq b$, where $M(a,b)$, $A(a,b)$ and
$C(a,b)$ denote the Neuman-S\'{a}ndor, arithmetic, and
contra-harmonic means of $a$ and $b$, respectively.
\end{abstract}

\bigskip
\centerline{\bf 1. Introduction} \setcounter{section}{1}
\setcounter{equation}{0}
\bigskip

For $a,b>0$ with $a\neq b$ the Neuman-S\'{a}ndor mean $M(a,b)$ [1] is defined by
\begin{equation*}
M(a,b)=\frac{a-b}{2{{\arcsinh}}\left[(a-b)/(a+b)\right]},
\end{equation*}
where ${\arcsinh}(x)=\log(x+\sqrt{1+x^2})$ is the inverse hyperbolic
sine function.

Recently, the Neuman-S\'{a}ndor mean has been the subject intensive
research. In particular, many remarkable inequalities for the
Neuman-S\'{a}ndor mean $M(a,b)$ can be found in the literature
[1-4].

Let $A(a,b)=(a+b)/2$, $G(a,b)=\sqrt{ab}$,
$L(a,b)=(b-a)/(\log{b}-\log{a})$, $C(a,b)=(a^2+b^2)/(a+b)$,
$P(a,b)=(a-b)/(4\arctan{\sqrt{a/b}}-\pi)$ and
$T(a,b)=(a-b)/\left[2\arctan((a-b)/(a+b))\right]$ be the arithmetic,
geometric, logarithmic, contra-harmonic, first Seiffert and second
Seiffert means of $a$ and $b$, respectively. Then it is well known
that the inequalities
\begin{equation*}
G(a,b)<L(a,b)<P(a,b)<A(a,b)<T(a,b)<C(a,b)
\end{equation*}
hold for all $a,b>0$ with $a\neq b$.

Neuman and S\'{a}ndor [1, 2] established that
\begin{equation*}
A(a,b)<M(a,b)<T(a,b)
\end{equation*}
\begin{equation*}
P(a,b)M(a,b)<A^{2}(a,b)
\end{equation*}
\begin{equation*}
A(a,b)T(a,b)<M^{2}(a,b)<(A^{2}(a,b)+T^{2}(a,b))/2
\end{equation*}
for all $a,b>0$ with $a\neq b$.

Let $0<a,b<1/2$ with $a\neq b$, $a'=1-a$ and $b'=1-b$. Then the
following Ky Fan inequalities
\begin{equation*}
\frac{G(a,b)}{G(a',b')}<\frac{L(a,b)}{L(a',b')}<\frac{P(a,b)}{P(a',b')}<\frac{A(a,b)}{A(a',b')}<\frac{M(a,b)}{M(a',b')}<\frac{T(a,b)}{T(a',b')}
\end{equation*}
were presented in [1].

Li et al. [3] proved that $L_{p_{0}}(a,b)<M(a,b)<L_{2}(a,b)$ for all
$a,b>0$ with $a\neq b$, where
$L_{p}(a,b)=[(b^{p+1}-a^{p+1})/((p+1)(b-a))]^{1/p}(p\neq -1, 0)$,
$L_{0}(a,b)=1/e(b^{b}/a^{a})^{1/(b-a)}$ and
$L_{-1}(a,b)=(b-a)/(\log{b}-\log{a})$ is the $p$-th generalized
logarithmic mean of $a$ and $b$, and $p_{0}=1.843\cdots$ is the
unique solution of the equation $(p+1)^{1/p}=2\log(1+\sqrt{2})$.
And, in [4] the author proved that the double inequality
\begin{equation}
\alpha C(a,b)+(1-\alpha) A(a,b)<M(a,b)<\beta C(a,b)+(1-\beta) A(a,b)
\end{equation}
holds for all $a,b>0$ with $a\neq b$ if and only if $\alpha \leq
\left(1-\log(\sqrt{2}+1)\right)/\log(\sqrt{2}+1)=0.1345\cdots$ and
$\beta \geq 1/6$, and the inequality
\begin{equation}
C^{\lambda}(a,b)A^{1-\lambda}(a,b)<M(a,b)<C^{\mu}(a,b)A^{1-\mu}(a,b)
\end{equation}
holds true for all $a,b>0$ with $a\neq b$ if $\mu\geq
\log\left((\sqrt{2}+2)/3\right)/\log{2}=0.1865\cdots$ and
$\lambda\leq 1/6$.

The main purpose of this paper is to give some refinements and
improvements for inequalities (1.1) and (1.2). Our main results are
the following Theorems 1.1 and 1.2.

\medskip
{\bf THEOREM 1.1.} The double inequality
\begin{equation*}
C^{\alpha}(a,b)A^{1-\alpha}(a,b)<M(a,b)<C^{\beta}(a,b)A^{1-\beta}(a,b)
\end{equation*}
holds for all $a,b>0$ with $a\neq b$ if and only if $\alpha\leq 1/6$
and $\beta \geq -\log(\log(1+\sqrt{2}))/\log{2}=0.1821\cdots$.

\medskip
{\bf THEOREM 1.2.} The double inequality
\begin{align*}
&[C(a,b)/6+5 A(a,b)/6]^{\lambda
}\left[C^{1/6}(a,b)A^{5/6}(a,b)\right]^{1-\lambda}<M(a,b)\\
&\qquad<[C(a,b)/6+5
A(a,b)/6]^{\mu}\left[C^{1/6}(a,b)A^{5/6}(a,b)\right]^{1-\mu}
\end{align*}
holds for all $a,b>0$ with $a\neq b$ if and only if $\lambda\leq
-[6\log(\log(1+\sqrt{2}))+\log{2}]/[6\log(7/6)-\log{2}]=0.27828\cdots$
and $\mu \geq 8/25$.

\bigskip
\bigskip
\centerline{\bf 2. Lemmas}
\setcounter{section}{2}\setcounter{equation}{0}
\bigskip

In order to prove our main results we need three Lemmas, which we
present in this section.

\medskip
{\bf LEMMA 2.1} (See [5, Theorem 1.25]). For $-\infty<a<b<\infty$,
let $f,g:[a,b]\rightarrow{\mathbb{R}}$ be continuous on $[a,b]$, and
be differentiable on $(a,b)$, let $g'(x)\neq 0$ on $(a,b)$. If
$f^{\prime}(x)/g^{\prime}(x)$ is increasing (decreasing) on $(a,b)$,
then so are
$$\frac{f(x)-f(a)}{g(x)-g(a)}\ \ \mbox{and}\ \ \frac{f(x)-f(b)}{g(x)-g(b)}.$$
If $f^{\prime}(x)/g^{\prime}(x)$ is strictly monotone, then the
monotonicity in the conclusion is also strict.

\medskip
{\bf LEMMA 2.2} (See [6, Lemma 1.1]). Suppose that the power series
$f(x)=\sum\limits_{n=0}^{\infty}a_{n}x^{n}$ and
$g(x)=\sum\limits_{n=0}^{\infty}b_{n}x^{n}$ have the radius of
convergence $r>0$ and $b_{n}>0$ for all $n\in\{0,1,2,\cdots\}$. Let
$h(x)={f(x)}/{g(x)}$, then the following statements are true:

(1) If the sequence $\{a_{n}/b_{n}\}_{n=0}^{\infty}$ is (strictly)
increasing (decreasing), then $h(x)$ is also (strictly) increasing
(decreasing) on $(0,r)$;

(2) If the sequence $\{a_{n}/b_{n}\}$ is (strictly) increasing
(decreasing) for $0<n\leq n_{0}$ and (strictly) decreasing
(increasing) for $n>n_{0}$, then there exists $x_{0}\in(0,r)$ such
that $h(x)$ is (strictly) increasing (decreasing) on $(0,x_{0})$ and
(strictly) decreasing (increasing) on $(x_{0},r)$.

\medskip
{\bf LEMMA 2.3.} The function
\begin{equation}
h(t)=\frac{90t+52t\cosh(2t)-66\sinh(2t)+2t\cosh(4t)-3\sinh(4t)}{15t-20t\cosh(2t)+5t\cosh(4t)}
\end{equation}
is strictly decreasing on $(0,\log(1+\sqrt{2}))$, where
$\sinh(t)=(e^{t}-e^{-t})/2$ and $\cosh(t)=(e^{t}+e^{-t})/2$ are the
hyperbolic sine and cosine functions, respectively.

\medskip
{\bf Proof.} Let
\begin{equation}
h_{1}(t)=90t+52t\cosh(2t)-66\sinh(2t)+2t\cosh(4t)-3\sinh(4t),
\end{equation}
\begin{equation}
h_{2}(t)=15t-20t\cosh(2t)+5t\cosh(4t).
\end{equation}
Then making use of power series formulas we have
\begin{align}
h_{1}(t)=&90t+52t\sum_{n=0}^{\infty}\frac{(2t)^{2n}}{(2n)!}-66\sum_{n=0}^{\infty}\frac{(2t)^{2n+1}}{(2n+1)!}
+2t\sum_{n=0}^{\infty}\frac{(4t)^{2n}}{(2n)!}-3\sum_{n=0}^{\infty}\frac{(4t)^{2n+1}}{(2n+1)!}\nonumber\\
=&52t\sum_{n=2}^{\infty}\frac{(2t)^{2n}}{(2n)!}-66\sum_{n=2}^{\infty}\frac{(2t)^{2n+1}}{(2n+1)!}
+2t\sum_{n=2}^{\infty}\frac{(4t)^{2n}}{(2n)!}-3\sum_{n=2}^{\infty}\frac{(4t)^{2n+1}}{(2n+1)!}\nonumber\\
=&\sum_{n=0}^{\infty}\frac{[16+13n+(2n-1)2^{2n+2}]2^{2n+7}}{(2n+5)!}t^{2n+5}
\end{align}
and
\begin{align}
h_{2}(t)=&15t-20t\sum_{n=0}^{\infty}\frac{(2t)^{2n}}{(2n)!}+5t\sum_{n=0}^{\infty}\frac{(4t)^{2n}}{(2n)!}\nonumber\\
=&-20t\sum_{n=2}^{\infty}\frac{(2t)^{2n}}{(2n)!}+5t\sum_{n=2}^{\infty}\frac{(4t)^{2n}}{(2n)!}=\sum_{n=0}^{\infty}\frac{5
(2^{2n+2}-1)2^{2n+6}}{(2n+4)!}t^{2n+5}.
\end{align}

It follows from (2.1)-(2.5) that
\begin{equation}
h(t)=\frac{\sum\limits_{n=0}^{\infty}a_{n}t^{2n}}{\sum\limits_{n=0}^{\infty}b_{n}t^{2n}},
\end{equation}
where
\begin{equation}
a_{n}=\frac{[16+13n+(2n-1)2^{2n+2}]2^{2n+7}}{(2n+5)!},\quad
b_{n}=\frac{5(2^{2n+2}-1) 2^{2n+6}}{(2n+4)!}.
\end{equation}

Equation (2.7) leads to
\begin{equation}
\frac{a_{n+1}}{b_{n+1}}-\frac{a_{n}}{b_{n}}=-\frac{6c_{n}}{5(2n+5)(2n+7)(2^{2n+2}-1)(2^{2n+4}-1)},
\end{equation}
where
\begin{equation}
c_{n}=(30n^2+135n+110-4^{n+3})4^{n+1}+11.
\end{equation}

From (2.9) we get
\begin{equation}
c_{0}=195,\quad c_{1}=315,\quad c_{2}=-33525
\end{equation}
and
\begin{align}
c_{n}&<(30n^2+135n+110-64n^3)4^{n+1}+11\nonumber\\
&=\left[10n^2(3-n)+15n(9-n^2)+5(22-n^3)-34n^3\right]4^{n+1}+11\nonumber\\
&<-34n^3\cdot 4^{n+1}+11<0
\end{align}
for $n\geq 3$.

Equations (2.8) and (2.10) together with inequality (2.11) lead to
the conclusion that the sequence $\{a_{n}/b_{n}\}$ is strictly
decreasing for $0\leq n \leq 2$ and strictly increasing for $n\geq
3$. Then from Lemma 2.2(2) and (2.6) we clearly see that there
exists $t_{0}\in(0,\infty)$ such that $h(t)$ is strictly decreasing
on $(0,t_{0})$ and strictly increasing on $(t_{0},\infty)$.

Let $t^{*}=\log(1+\sqrt{2})$. Then simple computations lead to
\begin{equation}
\sinh(2t^{*})=2\sqrt{2},\cosh(2t^{*})=3, \sinh(4t^{*})=12\sqrt{2},
\cosh(4t^{*})=17.
\end{equation}

Differentiating (2.1) yields
\begin{align}
h'(t)=&\frac{90-80\cosh(2t)+104t\sinh(2t)-10\cosh(4t)+8t\sinh(4t)}{h_{2}(t)}\nonumber\\
&-\frac{15-20\cosh(2t)-40t\sin(2t)+5\cosh(4t)+20t\sinh(4t)}{{h_{2}(t)}^2}h_{1}(t).
\end{align}

From (2.2) and (2.3) together with (2.12) and (2.13) we get
\begin{equation}
h'(t^*)=\frac{-102\sqrt{2}{t^{*}}^2+93t^{*}+21\sqrt{2}}{5{t^{*}}^2}=-0.10035\cdots<0.
\end{equation}

From the piecewise monotonicity of $h(t)$ and inequality (2.14) we
clearly see that $t_{0}>t^{*}=\log(1+\sqrt{2})$, and the proof of
Lemma 2.3 is completed. $\Box$

\bigskip
\bigskip
\centerline{\bf 3. Proof of Theorems 1.1 and 1.2}
\setcounter{section}{3}\setcounter{equation}{0}
\bigskip

{\bf\em Proof of Theorem 1.1.} Since $M(a,b)$, $C(a,b)$ and $A(a,b)$
are symmetric and homogeneous of degree $1$. Without loss of
generality, we assume that $a>b$. Let $x=(a-b)/(a+b)$ and
$t={\arcsinh}(x)$. Then $x\in(0,1)$,  $t\in(0,\log(1+\sqrt{2}))$ and
\begin{equation}
\frac{\log{[M(a,b)]}-\log{[A(a,b)]}}{\log{[C(a,b)]}-\log{[A(a,b)]}}=
\frac{\log[{x}/{{\arcsinh}(x)}]}{\log(1+x^2)}=\frac{\log[{\sinh(t)}/{t}]}{2\log[\cosh(t)]}.
\end{equation}

Let $f_{1}(t)=\log[{\sinh(t)}/{t}]$, $f_{2}(t)=\log[\cosh(t)]$ and
\begin{equation}
f(t)=\frac{\log[{\sinh(t)}/{t}]}{\log[\cosh(t)]}.
\end{equation}
Then $f_{1}(0^+)=f_{2}(0)=0$, $f(t)=f_{1}(t)/f_{2}(t)$ and
\begin{align}
\frac{{f_{1}}'(t)}{{f_{2}}'(t)}=&\frac{t\cosh^{2}(t)-\sinh(t)\cosh(t)}{t\sinh^{2}(t)}=\frac{t[\cosh(2t)+1]-\sinh(2t)}{t[\cosh(2t)-1]}\nonumber\\
=&\frac{t\left(\sum\limits_{n=0}^{\infty}2^{2n}t^{2n}/(2n)!+1\right)-\sum\limits_{n=0}^{\infty}2^{2n+1}t^{2n+1}/(2n+1)!}
{t\sum\limits_{n=1}^{\infty}2^{2n}t^{2n}/(2n)!}\nonumber\\
=&\frac{\sum\limits_{n=1}^{\infty}2^{2n}t^{2n+1}/(2n)!-\sum\limits_{n=1}^{\infty}2^{2n+1}t^{2n+1}/(2n+1)!}
{t\sum\limits_{n=1}^{\infty}2^{2n}t^{2n}/(2n)!}
=\frac{\sum\limits_{n=0}^{\infty}A_{n}t^{2n}}{\sum\limits_{n=0}^{\infty}B_{n}t^{2n}},
\end{align}
where $A_{n}=2^{2n+2}(2n+1)/(2n+3)!$ and
$B_{n}={2^{2n+2}}/{(2n+2)!}$.

Note the $A_{n}/B_{n}=1-2/(2n+3)$ is strictly increasing for all
$n\geq 0$. Then from Lemma 2.2(1) and (3.3) we know that
${f_{1}}'(t)/{f_{2}}'(t)$ is strictly increasing on $(0,\infty)$.
Hence, $f(t)$ is strictly increasing on $(0,\log(1+\sqrt{2}))$
follows from  Lemma 2.1 and the monotonicity of
${f_{1}}'(t)/{f_{2}}'(t)$ together with $f(0^+)=f_{2}(0)=0$.
Moreover,
\begin{equation}
\lim\limits_{t\rightarrow 0}f(t)=\lim\limits_{t\rightarrow
0}\frac{{f_{1}}'(t)}{{f_{2}}'(t)}=\frac{A_{0}}{B_{0}}=\frac{1}{3},
\end{equation}
\begin{equation}
\lim\limits_{t\rightarrow
\log(1+\sqrt{2})}f(t)=-\frac{2\log(\log(1+\sqrt{2}))}{\log{2}}.
\end{equation}

Therefore, Theorem 1.1 follows easily from (3.1), (3.2), (3.4) and
(3.5) together with the monotonicity of $f(t)$. $\Box$

\bigskip
{\bf\em Proof of Theorem 1.2.} Since $M(a,b)$, $C(a,b)$ and $A(a,b)$
are symmetric and homogeneous of degree $1$. Without loss of
generality, we assume that $a>b$. Let $x=(a-b)/(a+b)$ and
$t={\arcsinh}(x)$. Then $x\in(0,1)$,  $t\in(0,\log(1+\sqrt{2}))$ and
\begin{align}
&\frac{\log{M(a,b)}-\log{\left[C^{1/6}(a,b)A^{5/6}(a,b)\right]}}{\log{\left[C(a,b)/6+5A(a,b)/6\right]}-\log{\left[C^{1/6}(a,b)A^{5/6}(a,b)\right]}}\nonumber\\
=&\frac{\log[{x}/{{\arcsinh}(x)}]-\log(1+x^2)^{1/6}}{\log(1+x^2/6)-\log(1+x^2)^{1/6}}
=\frac{\log[{\sinh(t)}/{t}]-[\log\cosh(t)]/3}{\log[1+\sinh^{2}(t)/6]-[\log\cosh(t)]/3}.
\end{align}

Let $g_{1}(t)=\log[{\sinh(t)}/{t}]-[\log\cosh(t)]/3$,
$g_{2}(t)=\log[1+\sinh^{2}(t)/6]-[\log\cosh(t)]/3$ and
\begin{equation}
g(t)=\frac{\log[{\sinh(t)}/{t}]-[\log\cosh(t)]/3}{\log[1+\sinh^{2}(t)/6]-[\log\cosh(t)]/3}.
\end{equation}
Then $g_{1}(0^+)=g_{2}(0)=0$, $g(t)=g_{1}(t)/g_{2}(t)$ and
\begin{equation*}
\frac{{g_{1}}'(t)}{{g_{2}}'(t)}=\frac{[6+\sinh^{2}(t)][3t\cosh^{2}(t)-3\cosh(t)\sinh(t)-t\sinh^{2}(t)]}
{t\sinh(t)[6\sinh(t)\cosh^{2}(t)-\sinh(t)(6+\sinh^{2}(t))]}.
\end{equation*}

Elementary computations lead to
\begin{align*}
&[6+\sinh^{2}(t)][3t\cosh^{2}(t)-3\cosh(t)\sinh(t)-t\sinh^{2}(t)]\nonumber\\
=&\frac{45}{4}t+\frac{13}{2}t\cosh(2t)-\frac{33}{4}\sinh(2t)+\frac{t}{4}\cosh(4t)-\frac{3}{8}\sinh(4t),
\end{align*}
\begin{align*}
&t\sinh(t)[6\sinh(t)\cosh^{2}(t)-\sinh(t)(6+\sinh^{2}(t))]\nonumber\\
=&\frac{15}{8}t-\frac{5}{2}t\cosh(2t)+\frac{5}{8}t\cosh(4t)
\end{align*}
and
\begin{equation}
\frac{{g_{1}}'(t)}{{g_{2}}'(t)}=h(t),
\end{equation}
where $h(t)$ is defined as in Lemma 2.3.

It follows from Lemmas 2.1 and 2.3 and (3.8) together with
$g_{1}(0^+)=g_{2}(0)=0$ that $g(t)$ is strictly decreasing on
$(0,\log(1+\sqrt{2}))$. Moreover,
\begin{equation}
\lim\limits_{t\rightarrow 0}g(t)=\frac{8}{25},
\end{equation}
\begin{equation}
\lim\limits_{t\rightarrow
\log(1+\sqrt{2})}g(t)=-\frac{6\log(\log(1+\sqrt{2}))+\log{2}}{6\log(7/6)-\log{2}}.
\end{equation}

Therefore, Theorem 1.2 follows easily from (3.6), (3.7), (3.9) and
(3.10) together with the monotonicity of $g(t)$. $\Box$

\medskip
\def\refname
\bigskip
{\small \hfill {\em Miaokun Wang}

\hfill {\em Department of Mathematics}

\hfill {\em Huzhou Teachers College}

\hfill {\em Huzhou 313000}

\hfill {\em China}

\hfill {\em e-mail:} wmk000@126.com

\bigskip
{\small \hfill {\em Yuming Chu}

\hfill {\em Department of Mathematics}

\hfill {\em Huzhou Teachers College}

\hfill {\em Huzhou 313000}

\hfill {\em China}

\hfill {\em e-mail:} chuyuming@hutc.zj.cn

\bigskip
{\small \hfill {\em Baoyu Liu}

\hfill {\em School of Science}

\hfill {\em Hangzhou Dianzi University}

\hfill {\em Hangzhou 310018}

\hfill {\em China}

\hfill {\em e-mail:} 627847649@qq.com

{\hfil\bf\normalsize REFERENCES}
\begin{thebibliography}{ab}
\bibitem[1]{1} E. NEUMAN AND J. S\'{A}NDOR, \emph{On the Schwab-Borchardt mean}, Math. Pannon. {\bf 14}, 2(2003), 253-266.


\bibitem[2]{2} E. NEUMAN AND J. S\'{A}NDOR, \emph{On the Schwab-Borchardt mean II}, Math. Pannon. {\bf 17}, 1(2006), 49-59.

\bibitem[3]{3} Y.M. LI, B.Y. LONG AND Y.M. CHU, \emph{Sharp bounds for the Neuman-S\'{a}ndor mean in terms of generalized logarithmic mean}, J. Math. Inequal. {\bf 6}, 4(2012), 567-577.

\bibitem[4]{4} E. NEUMAN, \emph{A note on certain bivariate mean}, J. Math. Inequal. {\bf 6}, 4(2012), 637-643.

\bibitem[5]{5} G. D. ANDERSON, M. K. VAMANAMURTHY, M. VUORINEN, \emph{Conformal
Invariants, Inequalities, and Quasiconformal Maps}, John Wiley \&
Sons, New York, 1997.

\bibitem[6]{6} S. SIMI\'{C}, M. VUORINEN, \emph{Landen inequalities for zero-balanced
hypergeometric functions}, Abstr. Appl. Anal., 2012, Art. ID 932061,
11 pages.

\end{thebibliography}
\end{document}